\documentclass[10pt,a4paper,reqno]{amsart}
\usepackage{latexsym,amsmath,amssymb,paralist,graphicx}
\usepackage[bookmarks,ps2pdf,
pagebackref]{hyperref}

\author{Folkmar Bornemann and Christian Rasch}
\title[Finite-Element Discretization of Static Hamilton-Jacobi Equations]{Finite-Element
Discretization of Static Hamilton-Jacobi Equations Based on a Local Variational Principle}

\renewcommand{\MR}[1]{\MRhref{#1}{\relax\ifhmode\unskip\space\fi MR #1.}}
\renewcommand{\MRhref}[2]{%
  \href{http://www.ams.org/mathscinet-getitem?mr=#1}{#2}
}

\newcommand{\R}{\mathbb{R}}

\renewcommand{\leq}{\leqslant}
\renewcommand{\geq}{\geqslant}

\newcommand{\diam}{\operatorname{diam}}
\newcommand{\norm}[1]{\left\|#1\right\|}

\newcommand{\ska}[2]{\left\langle #1, #2\right\rangle}
\newcommand{\tol}{\text{\rm tol}}

\newtheorem{theorem}{Theorem}
\newtheorem{proposition}[theorem]{Lemma}

\begin{document}

\begin{abstract}
We propose a linear finite-element discretization of Dirichlet problems for static Hamilton-Jacobi equations
on unstructured triangulations. The discretization is based on simplified
localized Dirichlet problems that are solved by a local variational principle. It generalizes several
approaches known in the literature and allows for a simple and transparent convergence theory.
In this paper the resulting system of nonlinear equations is solved by an
{\em adaptive} Gauss--Seidel iteration that is easily implemented and quite effective as a couple of numerical experiments
show.
\end{abstract}
\address{Center of Mathematics, Technical University of Munich, 80290 Munich, Germany}
\email{\{bornemann,rasch\}@ma.tum.de}

\maketitle

\section{Introduction}

With the advent \cite{0659.65132} and success of level set methods and its many applications  \cite{Osher,Sethian}  to areas
ranging from computational physics to computer vision there has been considerable interest in
numerical methods for solving Hamilton--Jacobi equations, dynamic and static. For problems with complex
geometries, or for problems on manifolds, there is a demand for methods that work on unstructured meshes
such as triangulations.

 Three main directions of constructing discretizations on unstructured meshes can be found in the literature.
 First, there are methods that lift ideas of finite-difference upwinding and Godunov schemes from hyperbolic
 conservation laws to Hamilton--Jacobi equations (whose solutions are, at least in 1D, {\em integrals} of
 solutions to conservation laws), see, e.g., \cite{BS}. Second, there are finite-element methods that are
 based on a weak formulation of the semi-linear second order equation obtained by adding a small amount of factual
 viscosity, see, e.g., \cite{Xiang}.
 And third, there are methods that utilize the connection
 of Hamilton--Jacobi equations (via Bellman's  principle) to optimal control problems, see, e.g., \cite{SETHVLADSIAM}.

In this paper we propose a discretization
that bears similarities with the last approach. We implicitly
construct a linear finite-element solution by requiring that it solves {\em locally} a {\em simplified}
equation with (local) boundary conditions given by the finite-element function itself. The simplified local equation
is then solved by a \emph{local variational principle}, the {\em Hopf--Lax formula}.

This simple discretization is interesting in various respects. First, we will show that it generalizes quite a few
approaches known in the literature. Second, it allows for an extremely simple, self-contained convergence theory.
In fact, the only results of the general theory
that we rely on are a uniqueness theorem for viscosity solutions and the theorem of Arzel\`{a}--Ascoli. Existence
of viscosity solutions will be shown \emph{in passing} by the convergence of the finite-element discretization.

By construction the discretization inherits structural properties of the viscosity solution of the Hamilton--Jacobi
equation such as a comparison principle. For each property we will carefully trace
the specific assumptions on the Hamiltonian and the boundary
data that are needed for proofs in the continuous and the discrete case. In particular, it is known \cite{0890.49011,0497.35001}
that the existence of a viscosity solution of the Dirichlet problem necessitates a {\em compatibility condition} on the
boundary data, which is basically a restrictive Lipschitz bound. However, this necessary condition gets barely any mention in
the literature on numerical methods, even in the formulation of convergence results, e.g.,  \cite[Thm.~7.7]{SETHVLADSIAM}.

Moreover, we propose in this paper a likewise
simple iterative method for solving the resulting nonlinear system of
equations, namely the \emph{adaptive Gauss--Seidel iteration}, which is easily implemented and,
at least experimentally, quite effective.

The paper is organized as follows. In {\S}\ref{sect:visc} we recall the concept, existence, and uniqueness of viscosity solutions of
Dirichlet problems for certain Hamilton--Jacobi equations. In {\S}\ref{sect:rho} we introduce the support function of the zero-level
set of the Hamiltonian which plays a major role in the definition of the discretization. In {\S}\ref{sect:fem} we define the
linear finite-element solution based on a local variational principle. The existence, uniqueness, and uniform Lipschitz
continuity of the finite-element solutions are subject of {\S}\ref{sect:disc}. A suitable concept of {\em consistency} is introduced in
{\S}\ref{sect:conv} and the convergence of the discrete solutions is proved. In {\S}\ref{sect:eikonal} we apply the finite-element discretization to a class
of generalized eikonal equations in 2D and obtain, by a simple geometric argument, a closed formula for the local discrete
equation. In {\S}\ref{sect:gs} we discuss the adaptive Gauss--Seidel iteration that we propose for an easily implemented and quite effective solution of
the nonlinear system of equations. Finally, in {\S}\ref{sect:exp} we study two numerical experiments and compare the proposed method
to the ordered upwind method (OUM) recently published by Sethian and Vladimirsky \cite{SETHVLADSIAM}.

\medskip\bigskip\section{Existence and Uniqueness of Viscosity Solutions}\label{sect:visc}

In this section we shortly review the existence and uniqueness theory for the Dirichlet problem of a Hamilton--Jacobi equation,
\begin{equation}\label{eq.Dirich}
H(x,Du(x))=0,\quad x \in \Omega,\qquad u|_{\partial\Omega} = g,
\end{equation}
where throughout the paper we will assume that $\Omega \subset \R^d$ is a bounded Lipschitz domain.
For convex Hamiltonians $H$ a sufficiently general set of assumptions is (see \cite[{\S}5.3]{0497.35001}):
\begin{itemize}
\smallskip\item[(H1)] \emph{(Continuity)} $H \in C(\overline\Omega\times \R^d)$.
\smallskip\item[(H2)] \emph{(Convexity)} $p\mapsto H(x,p)$ is convex for all $x\in\overline\Omega$.
\smallskip\item[(H3)] \emph{(Coercivity)} $H(x,p)\rightarrow \infty$ as $\norm{p}\rightarrow \infty$, uniformly in $x\in\overline{\Omega}$. Equivalently,
by assumptions (H1) and (H2),
there are positive constants $\alpha, \beta$ with
\[
H(x,p) \geq \alpha \norm{p} - \beta,\qquad x \in \overline\Omega,\,p \in \R^d.
\]
\item[(H4)] \emph{(Compatibility of the Hamiltonian)} $H(x,0)\leq 0$ for all $x\in\overline\Omega$.
\end{itemize}
Existence of a solution to (\ref{eq.Dirich}) requires a further condition on the
boundary data:
\begin{itemize}
\smallskip\item[(H5)] \emph{(Compatibility of Dirichlet data)} $g(x) - g(y) \leq \delta(x,y)$ for all $x,y \in \partial\Omega$.
\end{itemize}
\smallskip
Here,
 $\delta$ denotes the {\em optical distance} defined, under the assumptions (H1)--(H4), by
\begin{multline}\label{eqn:DELTA}
\delta(x,y)=\inf \bigg\{ \int_0^1 \rho\big(\xi(t),-\xi'(t)\big) \, dt \, : \ \xi\in C^{0,1}([0,1],\overline\Omega),
\; \xi(0)=x, \ \xi(1)=y \bigg\}\\*[1mm]
\text{where}\qquad\rho(x,q)=\max_{H(x,p)=0}\ska{p}{q}.
\end{multline}
In fact, $\delta$ qualifies
as a distance by the fairly obvious properties
\begin{equation}\label{eq.dist}
\delta(x,x)=0,\qquad 0 \leq \delta(x,z)\leq \delta(x,y)+\delta(y,z),\qquad x,y,z \in \overline\Omega.
\end{equation}
If $H$ is symmetric with respect to $p$, that is, $H(x,p)=H(x,-p)$,
then $\delta$ defines a pseudometric
on $\overline{\Omega}$.

Let us recall the concept of viscosity solutions \cite{0543.35011} for the first order equation
\begin{equation}\label{eqn:HJ}
H(x,Du(x))=0,\qquad x \in \Omega.
\end{equation}
A function $u \in C^{0,1}(\overline\Omega)$
is a {\em viscosity subsolution (supersolution)} of (\ref{eqn:HJ}) if
 all
\mbox{$v\in C^\infty_0(\Omega)$} with $u-v$ attaining a local maximum (minimum) at some $x_0 \in \Omega$ yield
\[
H(x_0,Dv(x_0))\leq 0 \quad (\geq 0).
\]
Now, a {\em viscosity solution} is simultaneously a viscosity sub- and supersolution. Note that by
Rademacher's theorem on the differentiability of Lipschitz continuous functions a viscosity solution
satisfies (\ref{eqn:HJ}) pointwise almost everywhere.
\begin{theorem}[P.-L. Lions \protect{\cite[Thm.~5.3]{0497.35001}}]\label{thm:exi}
Assume (H1)--(H4). The Dirichlet problem (\ref{eq.Dirich}) has a viscosity solution $u$ if and only if the
boundary condition satisfies the compatibility condition (H5). A specific viscosity solution is then given by the Hopf--Lax formula
\begin{equation}\label{eqn:HopfLax}
u(x)=\inf_{y\in\partial\Omega} \Big( g(y) + \delta(x,y) \Big).
\end{equation}
\end{theorem}

While this theorem will only serve as a motivation for the finite-element discretization in {\S}\ref{sect:disc},
we will obtain the existence of viscosity solutions under somewhat more restrictive assumptions as a spin-off of the convergence result,
Theorem~\ref{thm:conv}.

Uniqueness requires a compatibility condition on the Hamiltonian that is slightly stronger than (H4):
\begin{itemize}
\smallskip\item[(H4$^\prime$)] $H(x,0)< 0$ for all $x\in\Omega$.
\end{itemize}
In fact, uniqueness of the viscosity solution given by (\ref{eqn:HopfLax}) is then
a simple corollary of the following {\em comparison principle}.
\begin{theorem}[H. Ishii \protect{\cite{0644.35017}}]\label{thm:comp}
Assume (H1)--(H3) and (H4$^\prime$). Let $u,v$ be viscosity sub- and supersolutions of (\ref{eqn:HJ}), respectively.
If $u\leq v$ on $\partial\Omega$ then $u\leq v$ on $\overline{\Omega}$.
\end{theorem}

\medskip\bigskip\section{The Support Function of the Zero-Level Set}\label{sect:rho}

In the definition (\ref{eqn:DELTA}) of the optical distance $\delta$ the {\em support function} (see \cite[p.~28]{0932.90001})
\begin{equation}\label{eqn:rho}
\rho(x,q) = \max_{H(x,p)=0} \langle p,q\rangle = \sup_{H(x,p)\leq 0} \langle p,q\rangle,\qquad x \in \overline\Omega,\; q \in \R^d,
\end{equation}
of the zero-level set of $H$ made an appearance.
It is a well-defined real-valued function, since by (H3) the zero-level set is
compact and by (H4) non-empty. The second equality in (\ref{eqn:rho}) follows from the convexity (H2).

Since the discretization that we propose in the next section will be based on this support function $\rho$
we collect its most important properties.

\begin{proposition}\label{lem:rho1}
Assume (H1)--(H4). Then $\rho : \overline\Omega\times \R^d \to \R$ is upper semicontinuous in the first argument,
positively homogeneous convex in the second,
that is,
\[
\rho(x,q_1+q_2)\leq \rho(x,q_1) + \rho(x,q_2),\qquad\rho(x,tq) = t\rho(x,q),\qquad t\geq 0,
\]
for $x\in \overline\Omega$ and $q,q_1,q_2 \in \R^d$. Let $\rho^*=\beta/\alpha$ with $\alpha$, $\beta$ from (H4). Then
\[
0 \leq \rho(x,q) \leq \rho^* \norm{q},\qquad x \in \overline\Omega,\; q \in \R^d.
\]
Assume additionally (H4$^\prime$). Then $\rho|_{\Omega\times\R^d}$ is continuous and
\[
\rho(x,q) > 0, \qquad x \in \Omega,\; 0 \neq q \in \R^d.
\]
\end{proposition}

\begin{proof} Being defined as the pointwise supremum of linear functions the function $q \mapsto \rho(x,q)$ is a convex,
positively homogeneous and, by (H4), nonnegative function for fixed $x \in \overline\Omega$. Assumption
(H3) yields for $H(x,p)\leq 0$ the bound $\norm{p} \leq \beta/\alpha = \rho^*$, which readily implies the upper bound on $\rho$.

To prove the upper semicontinuity let $x_n \to x_0 \in\overline \Omega$ and $q \in \R^d$. We
extract a subsequence $x_{n'}$ such that
\[
\rho(x_{n'},q) = \langle p_{n'},q\rangle \to \limsup_{n\to\infty}\rho(x_n,q),
\]
where $p_{n'}$ is a maximizing argument with $H(x_{n'},p_{n'})=0$. Because of the bound $\norm{p_{n'}}\leq \rho^*$ we can
assume without loss of generality that $p_{n'} \to p_0$. We obtain $H(x_0,p_0)=0$ and therefore
\[
\limsup_{n\to\infty} \rho(x_n,q) = \langle p_0,q\rangle \leq \rho(x_0,q).
\]

From now on, we assume (H4$^\prime$). Let $x\in\Omega$.
Since $H(x,0)<0$ there is $\delta>0$ such that
$H(x,p)\leq 0$ for $\norm{p}\leq \delta$. Thus, for $q \neq 0$
\[
\rho(x,q)\geq \max_{\norm{p}\leq \delta} \ska p q = \delta\norm{q} > 0.
\]

Finally, to prove the lower semicontinuity let $x_n \to x_0 \in \Omega$ and $q \in \R^d$.
There is a maximizing $p_0 \in \R^d$ with $\rho(x_0,q) = \langle p_0,q\rangle$ and $H(x_0,p_0)=0$. We extract a subsequence such that
$\rho(x_{n'},q)  \to \liminf_{n\to\infty}\rho(x_n,q)$ and, below,
construct a sequence $p_{n'} \to p_0$ with $H(x_{n'},p_{n'})\leq 0$. With it in hand we conclude
\[
\liminf_{n\to\infty}\rho(x_n,q) = \lim_{n'\to\infty}\rho(x_{n'},q) \geq \lim_{n'\to\infty} \langle p_{n'},q\rangle = \langle p_0,q\rangle = \rho(x_0,q).
\]
There is no loss of generality in assuming that either always $H(x_{n'},p_0)\leq 0$ or always $H(x_{n'},p_0)>0$.
In the first case we simply take $p_{n'}=p_0$. In the second case,
since $H(x_{n'},0)<0$, there is a $\lambda_{n'} \in (0,1)$ with $H(x_{n'},\lambda_{n'}p_0)= 0$ and we put $p_{n'}=\lambda_{n'}p_0$. We can assume that
$\lambda_{n'} \to \lambda_0 \in [0,1]$. Taking limits in
\[
0 = H(x_{n'},\lambda_{n'}p_0) \leq (1-\lambda_{n'}) H(x_{n'},0) + \lambda_{n'} H(x_{n'},p_0)
\]
yields $0 \leq (1-\lambda_0) H(x_0,0)$ which, by (H4$'$), implies $\lambda_0=1$ and $p_{n'}\to p_0$.
\end{proof}

Note that if (H4$'$) holds and $H$ is symmetric with respect to $p$,  then $q\mapsto\rho(x,q)$ defines
a norm on $\R^d$ for all $x\in\Omega$.

%

\begin{proposition}\label{lem:rho2}
Assume (H1)--(H4) and that the segment joining the points $y,z \in \overline\Omega$ belongs to $\overline\Omega$.
Let $\rho_*\geq 0$ be a constant such that $\rho(x,q) \geq \rho_* \|q\|$, $x \in \Omega$ and $q \in \R^d$.
Then, with the constant $\rho^*$ defined in Lemma~\ref{lem:rho1},
\[
\rho_* \norm{y-z} \leq \delta(y,z) \leq \rho^* \|y-z\|.
\]
If $H(x,p)$ does not depend on $x$, then $\rho(x,q)=\rho(q)$ does not depend on $x$ either and
\[
\delta(y,z) = \rho(y-z).
\]
\end{proposition}
\begin{proof}
The optical distance $\delta(y,z)$ is bounded by the expressions
\begin{multline*}
\rho_0\cdot \inf \bigg\{ \int_0^1 \norm{\xi'(t)} \, dt \, : \ \xi\in C^{0,1}([0,1],\overline\Omega)
\textrm{ such that }\xi(0)=y, \ \xi(1)=z \bigg\},
\end{multline*}
with $\rho_0=\rho_*$ for the lower bound and $\rho_0=\rho^*$ for the upper bound. The infimum is nothing but the minimal length
of a path joining the point $y$ and $z$ \emph{within} $\overline\Omega$. By the assumption on $y$ and $z$ this minimum is realized by the
segment joining them.

If $H$, and hence $\rho$, does not depend on $x \in \overline \Omega$, we obtain by Jensen's inequality for $\xi\in C^{0,1}([0,1],\overline\Omega)$
with $\xi(0)=y$ and $\xi(1)=z$ that
\[
\int_0^1 \rho(-\xi'(t))\,dt \geq \rho\left(-\int_0^1 \xi'(t)\,dt\right)= \rho(y-z).
\]
The lower bound is attained for the segment joining $y$ and $z$ yielding the assertion $\delta(y,z) = \rho(y-z)$
(see also \cite[Remark~5.7]{0497.35001}).
\end{proof}

\smallskip

\paragraph{\emph{Example.}} An important class\footnote{Which essentially covers the general case as we will see in Footnote~\ref{foot:rem} in {\S}\ref{sect:conv}.} of Hamiltonians satisfying (H1)--(H3) and (H4$'$) is given by
\[
H(x,p) = F(x,p) -1
\]
where $F \in C(\overline\Omega\times\R^d)$ is assumed to be positively homogeneous convex in $p$ with the bounds
\begin{equation}\label{eq:Fstar}
0 < F_* \leq F(x,p) \leq F^*,\qquad x\in \overline\Omega,\;\|p\|=1.
\end{equation}
Duality theory of nonnegative positively homogeneous convex function (\emph{gauges}) \cite[{\S}15]{0932.90001}
teaches that the support function $\rho$ of the zero-level set of $H$ is the
{\em polar} of $F$, that is,
\[
\rho(x,q) = \max_{p\neq 0} \frac{\langle p,q\rangle}{F(x,p)},\qquad F(x,p) = \max_{q\neq 0} \frac{\langle p,q\rangle}{\rho(x,q)}.
\]
Hence, for the (particular) Hamilton--Jacobi--Bellman equation \cite[Eq.~(22)]{SETHVLADSIAM} with
\begin{equation}\label{eq:HJB}
H(x,p) = \max_{\|q\|=1} \langle p,-q\rangle f(x,q) - 1,
\end{equation}
where $f$ is continuous with bounds $0< f_* \leq f(x,q) \leq f^*$, $x \in \overline\Omega$, $\|q\|=1$, we immediately read off that
\begin{equation}\label{eq:rhoHJB}
\rho(x,q) = \frac{\|q\|}{f(x,-q/\|q\|)},\qquad x\in \overline\Omega,\;q \in \R^d.
\end{equation}

\medskip\bigskip\section{The Finite-Element Discretization}\label{sect:fem}

\paragraph{\protect\emph{Linear Finite Elements.}}
Let us shortly recall the notion of linear finite-elements. For a sequence $h \to 0$ we consider a family $\Sigma_h$ of
shape-regular simplicial triangulations of (the now polytopal domain) $\Omega \subset \R^d$. We denote
the diameter of a (closed) simplex $\sigma \in \Sigma_h$ by $h_1(\sigma)$ and the minimal height of a vertex in $\sigma$ by $h_0(\sigma)$. We assume
\[
h = \max_{\sigma \in \Sigma_h} h_1(\sigma)
\]
and measure the shape-regularity by a uniform bound
\[
1 \leq \frac{h_1(\sigma)}{h_0(\sigma)} \leq \theta,\qquad \sigma\in \Sigma_h,\; h\to 0,
\]
where we call $\theta$ the \emph{regularity constant} of the family of triangulations.

The space of linear finite elements on $\Sigma_h$, that is,
continuous functions that are affine if restricted to a simplex $\sigma \in \Sigma_h$, is denoted by $V_h$ and
\[
I_h : C(\overline\Omega) \to V_h
\]
is the corresponding nodal interpolation operator. We endow $V_h$ with the maximum norm, that is, convergence in $V_h$
is the uniform convergence of the finite-element functions.

The set of nodal points (vertices) of the triangulation $\Sigma_h$ that belong
to $\overline \Omega$, $\Omega$, $\partial\Omega$ are denoted by $\overline\Omega_h$, $\Omega_h$, $\partial\Omega_h$, respectively.
Note that a finite-element function $u_h \in V_h$ is uniquely determined by its nodal values, that is, the values $u_h(x_h)$ for all
$x_h \in\overline\Omega_h$.

\begin{figure}[tbp]
\centerline{\includegraphics[width=0.7\textwidth]{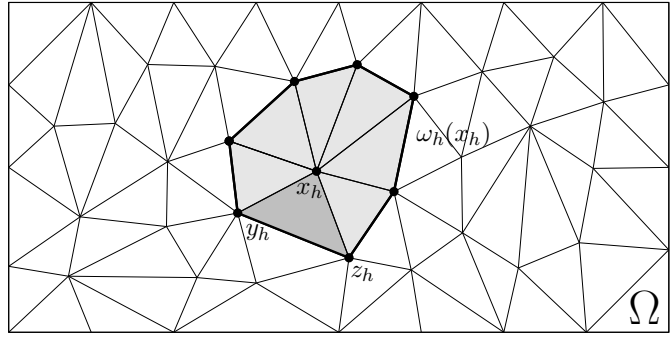}}
\caption{\protect\small The neighborhood $\omega_h(x_h)$ of $x_h \in \Omega_h$, that is, the collection of all simplices (like the one
shaded in dark) that have $x_h$ as a vertex.}
\label{fig:patch}
\end{figure}

For an interior nodal point $x_h \in \Omega_h$
we consider the simplicial neighborhood $\omega_h(x_h)$, that is,
the interior of the union of all simplices in $\Sigma_h$ that have $x_h$ as
a vertex (see Figure~\ref{fig:patch}).

\smallskip

\paragraph{\protect\emph{The Idea.}}
The finite-element discretization which we propose is motivated by the idea of {\em local solutions}:
\begin{quote}
At $x_h \in \Omega_h$ the finite-element solution $u_h \in V_h$ takes the value $u_h^*(x_h)$ of the {\em exact} viscosity solution
$u_h^* \in C^{0,1}(\overline{\omega_h(x_h)})$ that solves a
simplified Hamilton-Jacobi equation on $\omega_h(x_h)$ subject to the
 boundary conditions $u_h^*|_{\partial\omega_h(x_h)}=u_h|_{\partial\omega_h(x_h)}$.
\end{quote}
A good candidate for such a simplification of the Hamilton-Jacobi equation (\ref{eqn:HJ}) is obtained by freezing locally
the dependence of $H$ on its first variable. This way $u^*_h \in C^{0,1}(\overline{\omega_h(x_h)})$ is obtained as the viscosity solution of
the \emph{local Dirichlet problem}
\begin{equation}\label{eqn:HJlocal}
H(x_h,Du^*_h(x)) = 0 \;\;\text{on $\omega(x_h)$},\qquad u_h^*|_{\partial\omega_h(x_h)}=u_h|_{\partial\omega_h(x_h)}.
\end{equation}
This simplification is particularly suitable, because for $[y,z]\subset \overline{\omega_h(x_h)}$ the optical distance $\delta_{x_h}(y,z)$ of the local equation (\ref{eqn:HJlocal})
is, by Lemma~\ref{lem:rho2},
\[
\delta_{x_h}(y,z) = \rho(x_h,y-z),
\]
where $\rho(x_h,\cdot)$ is the support function of the zero-level set of the convex function $H(x_h,\cdot)$ as defined in (\ref{eqn:rho}).
Theorem~\ref{thm:exi} tells us that \emph{if} the local Dirichlet problem (\ref{eqn:HJlocal}) is solvable, then the value $u_h^*(x_h)$ is given
by the  {\em Hopf--Lax formula}
\[
u_h^*(x_h) = \min_{y \in \partial\omega_h(x_h)}\Big(u_h(y) + \rho(x_h,x_h-y) \Big),
\]
that is, by a simple {\em local variational principle}. We note that, under the assumptions (H1)--(H4), the Hopf--Lax formula is
well-defined \emph{independently} of
the compatibility of the boundary data of the local Dirichlet problem (\ref{eqn:HJlocal}).
Anyway, we base the finite-element discretization on this formula and the convergence will be proved later without
the interpretation of $u_h^*$ as a local solution.

\smallskip

\paragraph{\protect\emph{The Discretization}}
We define a function $\Lambda_h : V_h \to V_h$, called the {\em Hopf--Lax update} function, by
\[
(\Lambda_h u_h) (x_h) = \begin{cases}
\displaystyle\min_{y \in \partial\omega_h(x_h)}\Big(u_h(y) + \rho(x_h,x_h-y) \Big),\qquad & x_h \in \Omega_h,\\*[3mm]
u_h(x_h), & x_h \in \partial\Omega_h.
\end{cases}
\]
The finite-element solution $u_h \in V_h$ that discretizes the Dirichlet problem (\ref{eq.Dirich})
is now {\em implicitly} defined by the fixed-point equation
\begin{equation}\label{eq:fixedpoint}
u_h = \Lambda_h u_h,\qquad u_h|_{\partial\Omega_h} = g|_{\partial\Omega_h}.
\end{equation}
As in the continuous case, we call $u_h \in V_h$ a finite-element subsolution (supersolution) if $u_h \leq \Lambda_h u_h$ ($u_h \geq \Lambda_h u_h$).

We remark that the evaluation of $(\Lambda_h u_h)(x_h)$ at an interior nodal point $x_h \in \Omega_h$ can be calculated
by a finite collection
 of $(d-1)$-dimensional convex optimization problems. This follows from the representation (see Figure~\ref{fig:patch})
\begin{multline}\label{eqn:eval}
(\Lambda_h u_h) (x_h) = \min_{\sigma\in\Sigma_h : x_h \in \sigma} \min\Big\{u_h(y) + \rho(x_h,x_h-y) :\\
y \in \text{the ($d-1$)-dimensional face of $\sigma$
opposite to $x_h$}\Big\}
\end{multline}
and the observation that $u_h$ is affine, and hence $u_h + \rho(x_h,x_h-\cdot)$ convex, on $\sigma$.
In {\S}\ref{sect:eikonal} we will study an important class of examples for which these $(d-1)$-dimensional convex optimization
problems allow for a particularly simple solution. In general, however, one would have to use suitable iterative numerical methods
to solve them.

\smallskip
\paragraph{{\em Remark.}} For the particular Hamilton--Jacobi--Bellman equation (\ref{eq:HJB})
the finite-ele\-ment discretization (\ref{eq:fixedpoint}) is equivalent to various grid-based methods
that are obtained
from linear interpolation of the grid values and a direct local application of Bellman's dynamic programming principle.
See, e.g., \cite[Eq.~(2.3)]{0831.93028} and \cite[Eq.~(25)]{SETHVLADSIAM} as well as the references given therein.

\medskip\bigskip\section{Existence and Uniqueness of the Finite-Element Solution}\label{sect:disc}

The existence of a finite-element solution as implicitly defined by (\ref{eq:fixedpoint}) is based on two simple properties of the Hopf--Lax update function $\Lambda_h$.

\begin{proposition} Assume (H1)--(H4). Let $u_h, v_h \in V_h$.
\begin{enumerate}
\item $\Lambda_h$ is monotone, that is, $u_h \leq v_h$ implies $\Lambda_h u_h \leq \Lambda_h v_h$.\\*[-3mm]
\item $\Lambda_h$ is nonexpanding, that is, $\norm{\Lambda_h u_h - \Lambda_h v_h}_\infty \leq \norm{u_h - v_h}_\infty$.
\end{enumerate}
\end{proposition}
\begin{proof}
The first property is an immediate consequence of the definition of $\Lambda_h$. To prove the second, let the
maximum be attained at a nodal point $x_h \in \overline\Omega_h$, without loss of generality
$\norm{\Lambda_h u_h - \Lambda_h v_h}_\infty = (\Lambda_h u_h)(x_h) -(\Lambda_h v_h)(x_h)$.
If $x_h \in \partial\Omega_h$ there is nothing to show; so we can assume $x_h \in \Omega_h$. Let $y_* \in \partial\omega_h(x_h)$
be such that
\begin{equation}\label{eq:minHL}
(\Lambda_h v_h)(x_h) = v_h(y_*) + \rho(x_h,x_h-y_*).
\end{equation}
Hence
\begin{multline*}
(\Lambda_h u_h)(x_h) -(\Lambda_h v_h)(x_h) \\*[1mm]
\leq \Big(u_h(y_*) + \rho(x_h,x_h-y_*)\Big) - \Big(v_h(y_*) + \rho(x_h,x_h-y_*)\Big)
\leq \norm{u_h-v_h}_\infty,
\end{multline*}
which proves the assertion.
\end{proof}

\begin{theorem}\label{thm:exiFEM}
Assume (H1)--(H4) and $g : \partial \Omega \to \R$. Then the finite-element discretization (\ref{eq:fixedpoint})
has a solution $u_h \in V_h$. If $u^0_h \in V_h$ is such that $u^0_h|_{\partial\Omega_h} = g|_{\partial\Omega_h}$ and $\Lambda_h u_h^0 \geq u_h^0$,
then the fixed point iteration
\[
u_h^{n+1} = \Lambda_h u^n_h,\qquad n=0,1,2,\ldots,
\]
converges monotonously to a solution of (\ref{eq:fixedpoint}).
\end{theorem}
\begin{proof}
An initial iterate $u^0_h \in V_h$ with $\Lambda_h u_h^0 \geq u_h^0$ is given by
\[
u^0_h|_{\partial\Omega_h} = g|_{\partial\Omega_h},\qquad u_h|_{\Omega_h} = \min_{x \in \partial \Omega_h} g(x).
\]
Inductively the monotonicity of $\Lambda_h$ implies $u^{n+1}_h = \Lambda_h u^n_h \geq u_h^n$. Hence, the monotone convergence of the sequence
follows if we establish a uniform bound on the iterates. Since such a bound is trivial for the boundary we consider
for a given $x_h \in \Omega_h$ a shortest path $x_h = x_h^0,\ldots,x_h^m$ of nodal points that connects $x_h$ along edges of the
triangulation with the boundary: $x_h^m \in \partial \Omega_h$. There is a bound $L$ on the length of such a path
which depends on the triangulation but not on $x_h$.
With $\rho^*$ as defined in Lemma~\ref{lem:rho1} we get
\begin{equation}\label{eq:bound}
u_h^n(x_h) \leq \max_{y \in \partial\Omega_h} g(y) + \sum_{i=0}^{m-1} \rho(x_h^i,x_h^i-x_h^{i+1}) \leq
\max_{y \in \partial\Omega_h} g(y) + \rho^* L.
\end{equation}
Thus, $u_h^n \to u_h \in V_h$ for some $u_h \in V_h$, which by continuity must be a fixed point of $\Lambda_h$.
\end{proof}

Like in the continuous case, uniqueness of the finite-element solution requires the sharper condition (H4$'$) and is
a simple corollary of the following {\em discrete comparison principle}. Thus, the finite-element discretization is
a {\em monotone scheme}.

\begin{theorem}\label{thm:compFEM} Assume (H1)--(H3) and (H4$'$).
Let $u_h,v_h \in V_h$ be finite-element sub- and supersolutions, respectively.
If $u_h\leq v_h$ on $\partial\Omega_h$ then $u_h\leq v_h$ on $\overline \Omega$.
\end{theorem}
\begin{proof}
Let be $\Delta_h=u_h-v_h \in V_h$. Note that the maximum of $\Delta_h$ will be attained in a nodal point.
We will show that the existence of $x_h\in\Omega_h$
with $\Delta_h(x_h)=\max_{x\in\overline\Omega} \Delta_h(x)=\delta>0$ yields a contradiction.
To this end we choose such a maximizing $x_h$ with minimal value of $v_h(x_h)$.
With $y_*\in\partial\omega_h(x_h)$ as in (\ref{eq:minHL}) we get
\[
\delta=u_h(x_h) - v_h(x_h) \leq (\Lambda_h u_h)(x_h) - (\Lambda_h v_h) (x_h) \leq u_h(y_*) - v_h(y_*).
\]
On the boundary of the face that contains $y_*$ in its relative interior there is, by the maximality of $\delta$,
a point $x_h^* \in \Omega_h$  such that $\Delta_h(x_h^*) = \delta$ and $v_h(x_h^*) \leq v_h(y_*)$.
By Lemma~\ref{lem:rho1} we have $\rho(x_h,x_h- y_*)>0$ and obtain
\[
v_h(x_h^*) \leq v_h(y_*) = (\Lambda_h v_h)(x_h) - \rho(x_h,x_h- y_*) <  (\Lambda_h v_h)(x_h) \leq v_h(x_h)
\]
in contradiction to the minimality of $v_h(x_h)$.
\end{proof}

In the discrete case, up to now, we did not impose a compatibility condition on the boundary data such as (H5).
This will change in the discussion of a third important property of the finite-element solutions needed for convergence,
that is, uniform Lipschitz continuity. Based on the constant $\rho_*\geq 0$ of
Lemma~\ref{lem:rho2} and the regularity constant $\theta\geq 1$ of the family of triangulation we consider the condition
\begin{itemize}
\smallskip\item[(H5$'$)] $g(x) - g(y) \leq \dfrac{\rho_*}{\theta}\norm{x-y}$ for all $x,y \in \partial\Omega$.
\end{itemize}\smallskip
By Lemma~\ref{lem:rho2} this condition is actually stronger than (H5). Note that the {\em homogeneous} Dirichlet condition $g = 0$
always satisfies (H5$'$).

\begin{theorem}\label{thm:lipFEM} Assume (H1)--(H3), (H4$'$) and (H5$'$).  The unique finite-element
solution $u_h \in V_h$ of (\ref{eq:fixedpoint})
satisfies the uniform Lipschitz condition
\[
|u_h(x) - u_h(y)| \leq c_\Omega\, \theta d\cdot\rho^* \cdot \norm{x-y}, \qquad x,y \in \overline\Omega,
\]
and the uniform bound
\[
\| u_h \|_\infty \leq \max_{x \in \partial \Omega} |g(x)| + c_\Omega  \,\theta d\cdot \rho^*\cdot \diam(\Omega).
\]
Here, $\theta$ denotes the regularity constant of the family of triangulations, $\rho^*$ is the
constant defined in Lemma~\ref{lem:rho1} and $c_\Omega>0$ is a constant depending only on $\Omega$.
If $\Omega$ is convex, we can choose $c_\Omega = 1$.
\end{theorem}

\begin{proof} The uniform bound on $\norm{u_h}_\infty$ is a simple consequence of the Lipschitz condition.
The proof of the Lipschitz condition proceeds in three steps, imposing less and less restrictions on the possible choices of $x,y \in \overline\Omega$.

\paragraph{\emph{Step 1.}} For \emph{neighboring} nodal points $x_h,y_h \in \overline\Omega_h$ we prove
\[
|u_h(x_h) - u_h(y_h)| \leq \rho^* \cdot \norm{x_h-y_h}.
\]
Since $\rho^* \geq \rho_* \geq \rho_*/\theta$ this is, by (H5$'$), obviously true for $x_h,y_h \in \partial\Omega_h$.
If $x_h \in \Omega_h$ we have
$y_h \in \partial\omega_h(x_h)$ and hence
\[
u_h(x_h) = (\Lambda_h u_h)(x_h) \leq u_h(y_h) + \rho(x_h,x_h-y_h) \leq u_h(y_h) + \rho^* \norm{x_h-y_h}.
\]
If $y_h \in \Omega_h$ we can change the roles of $x_h$ and $y_h$ and the Lipschitz bound follows.

Assume on the other hand that $y_h \in \partial \Omega_h$. There is a minimizing $y_*\in\partial\omega_h(x_h)$ such that
\[
u_h(x_h) = (\Lambda_h u_h)(x_h) = u_h(y_*) + \rho(x_h,x_h-y_*) > u_h(y_*),
\]
where the last inequality follows from Lemma~\ref{lem:rho1}.
The boundary of the face that contains $y_*$ in its relative interior
has a point $x_h^1 \in \overline\Omega_h$ with
$u_h(x_h^1) \leq u_h(y_*) <u_h(x_h)$. By the definition of $\rho_*$ and $\theta$ we obtain
\[
\rho(x_h,x_h-y_*) \geq \rho_* \|x_h-y_*\| \geq \frac{\rho_*}{\theta} \|x_h-x_h^1\|.
\]
Continuing this construction we obtain a sequence
$x_h=x_h^0,x_h^1,\ldots,x_h^m$ of nodal points with strictly decreasing $u_h$-values that necessarily reaches the boundary at
some index $m$: $x_h^m \in \partial\Omega_h$. Thus, by construction and (H5$'$),
\begin{multline*}
u_h(x_h) \geq g(x_h^m) + \frac{\rho_*}{\theta}\sum_{i=0}^{m-1} \|x_h^i - x_h^{i+1}\| \\*[2mm]
\geq g(y_h) + \frac{\rho_*}{\theta}
\left(\sum_{i=0}^{m-1} \|x_h^i - x_h^{i+1}\| - \|x_h^m - y_h\|\right)\\*[2mm] \geq u_h(y_h) - \frac{\rho_*}{\theta} \norm{x_h-y_h}
\geq u_h(y_h) - \rho^* \norm{x_h-y_h},
\end{multline*}
which concludes the proof of Step 1.

\paragraph{\emph{Step 2.}} Let $\sigma \in \Sigma_h$ be a simplex of the triangulation. For $x,y \in\sigma$ we prove that
\[
|u_h(x) - u_h(y)| \leq \theta d\cdot\rho^*\cdot\norm{x-y}.
\]
By an affine transformation $\hat x \mapsto B \hat x + b$ we map the standard $d$-dimensional simplex
\[
\hat\sigma = \{\hat x\in \R^d_{\geq 0}: \hat x_1+\ldots +\hat x_d \leq 1\}
\]
onto $\sigma$.
The pullback of $u_h|_\sigma$ under the transformation will be denoted $\hat u$.
By Step~1 we can estimate the length of the (constant) gradient of $u_h|_\sigma$ by
\[
\|D u_h|_\sigma\| \leq \|B^{-1}\|\, \|D \hat u\| \leq \|B^{-1}\| \sqrt{d}\,\rho^*\cdot h_1(\sigma).
\]
Now, $\|B^{-1}\|$ is the largest ratio of the length of a segment in $\hat\sigma$ to the length of its image in $\sigma$.
Without loss of generality such a segment can be assumed to join a vertex with the opposite boundary face. Thus $\|B^{-1}\| \leq \sqrt{d}/h_0(\sigma)$ and,
by the shape-regularity assumption, that is, $h_1(\sigma)/h_0(\sigma) \leq \theta$, we get
\[
\|D u_h|_\sigma\| \leq \theta d\cdot \rho^*
\]
 and hence the assertion of Step~2.

\paragraph{\emph{Step 3.}} For $x,y \in \overline \Omega$ there is a Lipschitz path $\gamma \in C^{0,1}([0,1],\overline\Omega)$
joining $x$ and $y$ such that (see \cite[p.~304]{Alt02})
\[
\norm{\gamma'}_\infty \leq c_\Omega \|x-y\|.
\]
For convex $\Omega$ the path $\gamma$ can be chosen as the segment joining $x$ and $y$, which yields $c_\Omega=1$.
Now, let $0=t_0<t_1<\ldots<t_m=1$ be a subdivision of $[0,1]$ such that $\gamma(t_{i-1})$ and $\gamma(t_i)$ are elements of a common
simplex. By Step~2 we obtain
\begin{multline*}
|u_h(x)-u_h(y)| \leq \sum_{i=0}^{m-1} |u_h(\gamma(t_i))-u_h(\gamma(t_{i+1}))| \leq \theta d\cdot\rho^* \sum_{i=0}^{m-1}
\|\gamma(t_i)-\gamma(t_{i+1})\|\\*[2mm]
\leq c_\Omega\, \theta d\cdot\rho^*\cdot\|x-y\| \sum_{i=0}^{m-1} |t_i -t_{i+1}| =  c_\Omega\, \theta d\cdot\rho^*\cdot\|x-y\|,
\end{multline*}
which concludes the proof of the asserted Lipschitz bound.
\end{proof}

\medskip\bigskip\section{Convergence of the Finite-Element Discretization}\label{sect:conv}

The argument will be simplified if we consider a modified Hamiltonian $\tilde H$ for which the corresponding
Hamilton-Jacobi equation possesses the same viscosity solutions as the original one.

\begin{proposition}\label{prop:alt}
Assume (H1)--(H4). Let $x \in \bar\Omega$ and $p \in \R^d$. For the modified Hamiltonian
\[
\tilde H(x,p) = \max_{\|q\|=1}\left(\langle p,q\rangle - \rho(x,q)\right)
\]
we get that $\tilde H(x,p) \leq 0$ ($\tilde H(x,p) \geq 0$) implies $H(x,p) \leq 0$ ($H(x,p) \geq 0$).
\end{proposition}

\begin{proof}
First assume $H(x,p)>0$. There is a hyperplane
that separates $p$ strongly from the compact and convex level set $\{\tilde p : H(x,\tilde p) \leq 0\}$ (see \cite[Cor.~11.4.2]{0932.90001}). That is,
there is a vector $q \in \R^d$, $\norm{q}=1$,
with $\rho(x,q) < \langle p,q\rangle$.
Hence
\[
\tilde H(x,p) \geq \langle p,q\rangle - \rho(x,q) > 0.
\]
Now assume $H(x,p)<0$. There is $\epsilon > 0$ such that $H(x,p + \delta p)<0$ for $\norm{\delta p}\leq \epsilon$. Hence
\[
\langle p, q\rangle - \rho(x,q) \leq \langle p, q\rangle - \max_{\norm{\delta p}\leq \epsilon}\langle p + \delta p,q\rangle = - \epsilon \norm{q},
\]
Taking the supremum over all $q$ with $\norm{q}=1$ yields $\tilde H(x,p) \leq -\epsilon < 0$.
\end{proof}

In particular, each viscosity subsolution (supersolution) of the thus modified Hamilton--Jacobi equation
$\tilde H(x,Du(x)) = 0$, $x \in \Omega$,
is also a viscosity subsolution (supersolution) of the original one $H(x,Du(x)) = 0$, $x \in \Omega$.\footnote{\label{foot:rem}%
Under the additional assumption (H4$'$) the same holds true, since then $\rho(x,q)>0$ for $x\in\Omega$, if we consider
the modified Hamilton--Jacobi equation with the Hamiltonian
\[
\tilde H(x,p) = \max_{\|q\|=1}\frac{\ska{p}{q}}{\rho(x,q)} - 1.
\]
Hence, we see that the example at the end of {\S}\ref{sect:visc} in fact covers the general case.
}

Loosely speaking, in the framework of viscosity solutions the notion of \emph{consistency} of a discretization means that
a smooth
function is already a subsolution (supersolution) of the differential equation if it is a subsolution (supersolution) of
the discrete scheme. The precise statement is given in the next theorem.

\begin{theorem}\label{thm:consist} Assume (H1)--(H3) and (H4$'$).
Let $v\in C^\infty_0(\Omega)$, $x\in\Omega$, and $x_h\in\Omega_h$ be a sequence of nodal points
that converges to $x$ as $h\to 0$. Then
\begin{align*}
v(x_h) \leq (\Lambda_h I_h v) (x_h)\;\; \text{for all $h$} &\quad\Rightarrow\quad H(x,Dv(x))\leq 0,\\
v(x_h) \geq (\Lambda_h I_h v) (x_h)\;\; \text{for all $h$} &\quad\Rightarrow\quad H(x,Dv(x))\geq 0,
\end{align*}
where $I_h : C(\overline\Omega) \to V_h$ denotes the nodal interpolation operator.
\end{theorem}
\begin{proof} Since $v$ is smooth we can approximate the directional derivatives of $v$ in $x_h$ by first order differences
as follows
\begin{equation}\label{eq:interp}
\frac{v_h(x_h)-(I_hv)(y)}{\|x_h-y\|} = \left\langle Dv(x_h), \frac{x_h-y}{\norm{x_h-y}} \right\rangle  +
O(h),\qquad y \in \partial\omega_h(x_h).
\end{equation}
Now, let $v(x_h) \leq (\Lambda_h I_h v) (x_h)$ for all $h$ of the sequence, that is,
\[
v(x_h) - (I_h v)(y) - \rho(x_h,x_h-y) \leq 0,\qquad y \in \partial\omega_h(x_h).
\]
After division by $\|x_h-y\|$ we get, by (\ref{eq:interp}), a constant $c>0$ such that
\[
\langle Dv(x_h), q \rangle -  \rho(x_h,q) \leq c h,\qquad \|q\|=1.
\]
If we pass to the limit $h\to 0$ (note the continuity of $\rho$ at
$x \in \Omega$ as stated in Lemma~\ref{lem:rho1}) and take thereafter the maximum over all $\|q\|=1$, we obtain
\[
\tilde H(x,Dv(x)) = \max_{\|q\|=1}(\langle Dv(x), q \rangle -  \rho(x,q)) \leq 0.
\]
From Lemma~\ref{prop:alt} we infer the assertion $H(x,Dv(x)) \leq 0$.

On the other hand, let $v(x_h) \geq (\Lambda_h I_h v) (x_h)$ for all $h$ of the sequence, that is,
\[
v(x_h) - (I_h v)(y_h) - \rho(x_h,x_h-y_h) \geq 0
\]
for some $y_h \in \partial\omega_h(x_h)$. After devision by $\|x_h-y_h\|$ we get,
by (\ref{eq:interp}), a constant $c>0$ such that
\[
\langle Dv(x_h), q_h \rangle -  \rho(x_h,q_h) \geq -c h,\qquad q_h = (x_h-y_h)/\|x_h-y_h\|.
\]
By compactness, we can assume that $q_h \to q_*$ with $\|q_*\|=1$. Passing to the limit $h\to 0$ we thus obtain
\[
\tilde H(x,Dv(x)) \geq \langle Dv(x), q_* \rangle -  \rho(x,q_*) \geq 0,
\]
from which we infer the assertion $H(x,Dv(x))\geq 0$ by Lemma~\ref{prop:alt}.
\end{proof}

Now we have all the tools in hand to prove the \emph{convergence} of the finite-element discretization.

\begin{theorem}\label{thm:conv} Assume (H1)--(H3), (H4$'$), and (H5$'$). Then, as $h \to 0$,
the sequence of unique finite-element solutions $u_h \in V_h$ defined by
\[
u_h = \Lambda_h u_h,\qquad u_h|_{\partial\Omega_h} = g|_{\partial\Omega_h},
\]
converges uniformly to the unique viscosity solution $u$ of the Dirichlet problem
\[
H(x,Du(x)) = 0,\qquad u|_{\partial \Omega} = g.
\]
\end{theorem}
\begin{proof} Theorems~\ref{thm:exiFEM} and \ref{thm:compFEM} show the existence and uniqueness of the
finite-element solutions $u_h \in V_h$. Theorem~\ref{thm:lipFEM} shows that $u_h \in V_h$ is a uniform
bounded sequence of  uniform Lipschitz continuous functions. By
the theorem of Arzel\`{a}--Ascoli there is a subsequence $(u_{h'})$ that converges uniformly to
a function $u \in C^{0,1}(\overline\Omega)$. Because of  (H5$'$) and  $u_h|_{\partial\Omega_h} = g|_{\partial\Omega_h}$,
this limit satisfies the boundary condition $u|_{\partial \Omega} = g$.

To show that $u$ is a viscosity \emph{subsolution} of $H(x,Du(x)) = 0$ let $v\in C^\infty_0(\Omega)$ and
$x_0\in \Omega$ such that $u-v$ attains a local maximum in $x_0$. By adding a quadratic parabola to $v$ if necessary,
we may assume that it is in fact a strict local maximum (see \cite[p.~542]{Evans}). Extracting a further subsequence
of $h'$ if necessary, there is, by uniform convergence and
the monotonicity of the nodal interpolation operator $I_{h}: C(\overline \Omega) \to V_{h}$,
a sequence of nodal points $x_{h'} \in \Omega_{h'}$ such that $x_{h'}\to x_0$ and (see the argument given in \cite[p.~541]{Evans})
\[
(u_{h'}- v)(x_{h'}) \geq  (u_{h'}-I_{h'}v)(y),\qquad y\in\partial\omega_{h'}(x_{h'}).
\]
Now let $y_*\in\partial\omega_{h'}(x_{h'})$ be a minimizing argument such that
\[
(\Lambda_{h'} I_{h'}v) (x_{h'}) = (I_{h'} v)(y_*) + \rho(x_{h'},x_{h'}-y_*).
\]
Then it holds that
\begin{multline*}
(u_{h'} - v) (x_{h'})
 \geq  u_{h'}(y_*) + \rho(x_{h'},x_{h'}-y_*)  - (I_{h'}v)(y_*) - \rho(x_{h'},x_{h'}-y_*)\\
\geq (\Lambda_{h'} u_{h'}  - \Lambda_{h'} I_{h'} v) (x_{h'})
=(u_{h'}  - \Lambda_{h'} I_{h'}v_{h'} ) (x_{h'})
\end{multline*}
and thus
\[
v(x_{h'}) \leq (\Lambda_{h'} I_{h'}v ) (x_{h'}).
\]
The consistency of the discretization, stated in Theorem~\ref{thm:consist}, yields that
\[
H(x_0,Dv(x_0))\leq0,
\]
which concludes the proof that $u$ is a viscosity subsolution.

In the same way we prove that $u$ is a viscosity \emph{supersolution} of $H(x,Du(x)) = 0$. Therefore, $u$ is a viscosity solution,
which, by the comparison principle (Theorem~\ref{thm:comp}), is actually {\em unique}. Hence, there is exactly one limit point of
the sequence $u_h$, which thus has to converge uniformly to the just established viscosity solution $u$.
\end{proof}

\paragraph{\em Remark} Note that the only use that we have made so far of the existence
Theorem~\ref{thm:exi} was to motivate the local variational principle for the finite-element
discretization. In fact, our proof of the convergence result shows the existence of a viscosity solution {\em en route}
 --- under the somewhat more restrictive compatibility conditions (H4$'$) and (H5$'$), however.

\medskip\bigskip\section{The Hopf--Lax Update for Generalized Eikonal Equations in 2D}\label{sect:eikonal}

Let $\Omega \subset \R^2$ be a polygonal Lipschitz domain and $M:\overline\Omega \to \R^{2\times 2}$ be a continuous mapping into
the {\em symmetric positive definite} $2\times 2$-matrices. We denote the corresponding inner product by
$\langle p,q\rangle_{M(x)} = \langle M(x)p,q\rangle$, its subordinate norm by $\|p\|_{M(x)} = \langle p,p\rangle_{M(x)}^{1/2}$.

Now, we consider the Dirichlet problem for the \emph{generalized eikonal equation},
\[
\|Du\|_{M(x)} = 1\;\;\text{in $\Omega$},\qquad u|_{\partial\Omega} = g.
\]
Its Hamiltonian $H(x,p)=\|p\|_{M(x)}-1$ satisfies the assumptions (H1)--(H3) and (H4$'$). The support function of the zero-level
set is simply given by the norm that is dual to $\|\cdot\|_{M(x)}$, namely,
\[
\rho(x,q) = \max_{H(x,p)=0}\langle p,q \rangle = \max_{\|p\|_{M(x)}=1}\langle p,q \rangle = \|q\|_{M(x)^{-1}}.
\]
The Hopf--Lax update function becomes
\[
(\Lambda_h u_h)(x_h) = \min_{y \in \partial\omega_h(x_h)}\Big(u_h(y) + \|x_h-y\|_{M(x_h)^{-1}}\Big),\qquad x_h \in \Omega_h,\;u_h \in V_h.
\]
There is a simple procedure to evaluate $(\Lambda_h u_h)(x_h)$ at $x_h\in\Omega_h$. To this end let
$\sigma_1,\ldots,\sigma_m \in \Sigma_h$ be the triangles that have $x_h$ as a vertex and $J_i$ the (closed) edge of $\sigma_i$ opposite
to $x_h$. Then, as in (\ref{eqn:eval}),
\[
(\Lambda_h u_h)(x_h) = \min_{1\leq i \leq m} u_i\qquad\text{with}\qquad u_i = \min_{y \in J_i}\Big(u_h(y) + \|x_h-y\|_{M(x_h)^{-1}}\Big).
\]
Let us take one of the triangles, $\sigma_i$, (see Figure~\ref{fig:patch}) and call its vertices $x_h$, $y_h$, $z_h$, hence $J_i = [y_h,z_h]$.
In the case of the classic eikonal equation, that is, $M(x)\equiv I$, the update $u_i$ can be determined from an elementary geometric
argument.

\begin{proposition}\label{lem:update}
Let $\sigma \in \Sigma_h$ be the triangle with the vertices $x_h$, $y_h$, $z_h$ and $u_h \in V_h$. Denote the angles at $y_h$, $z_h$ by
$\alpha$, $\beta$, respectively. Defining
\[
\Delta = \frac{u_h(z_h)-u_h(y_h)}{\|z_h-y_h\|}
\]
and $\cos(\delta) = \Delta$ if $|\Delta|\leq 1$, we obtain
\begin{multline*}
u_i=\min_{y \in [y_h,z_h]} \Big(u_h(y) + \|x_h-y\| \Big) = u_h(y_h) + \min_{y \in [y_h,z_h]} \Big(\Delta\cdot\|y-y_h\| + \|x_h-y\| \Big)\\*[2mm]
=
\begin{cases}
u_h(y_h) + \|x_h-y_h\|,& \qquad \cos(\alpha) \leq \Delta,\\*[1mm]
u_h(y_h) + \cos(\delta-\alpha)\cdot\|x_h-y_h\|,& \qquad \alpha \leq \delta \leq \pi-\beta,\\*[1mm]
u_h(z_h) + \|x_h-z_h\|,& \qquad \Delta \leq \cos(\pi-\beta).
\end{cases}
\end{multline*}
\end{proposition}
\begin{figure}[tbp]
\centerline{\includegraphics[width=0.68\textwidth]{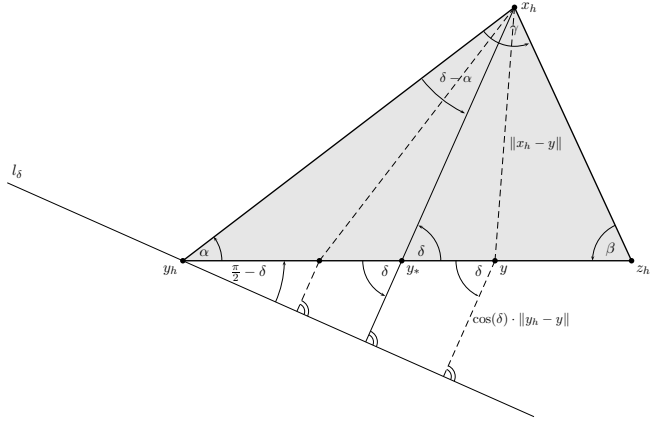}}
\caption{\protect\small Geometry of the minimization of $\cos(\delta)\|y-y_h\| + \|x_h-y\|$ for $y \in [y_h,z_h]$. Note that
for $\delta > \pi/2$ the segment through $y$ perpendicular to $l_\delta$ has \emph{negative} length $\cos(\delta)\cdot \|y-y_h\|$.}
\label{fig.proof}
\end{figure}%
\begin{proof}
For $|\Delta| \geq 1$ the assertion follows from a direct application of the triangle inequality; e.g., for $\Delta \geq  1$,
\[
\Delta\cdot\|y-y_h\| + \|x_h-y\| \geq \|y-y_h\| + \|x_h-y\| \geq \|x_h-y_h\|.
\]
Now, let $|\Delta|< 1$ so that $\cos(\delta) = \Delta$ defines a $\delta \in (0,\pi)$. A look at Figure~\ref{fig.proof} shows
that
\begin{equation}\label{eqn:simple}
\cos(\delta) \cdot\|y-y_h\| + \|x_h-y\|
\end{equation}
attains its minimum at the unique intersection $y_*$ of two straight lines: the first line running through $y_h$ and $z_h$, the
second line running through $x_h$ perpendicular to $l_\delta$. Here, $l_\delta$ is the straight line that encloses at $y_h$
with $[y_h,z_h]$ the angle $\pi/2-\delta$. We observe that the value of the minimum is simply $\cos(\delta-\alpha)\cdot \|x_h-y_h\|$.
A further look at Figure~\ref{fig.proof} teaches that $y_* \in [y_h,z_h]$ if and only if
\[
0 \leq \delta - \alpha \leq \gamma = \pi - \alpha - \beta,\qquad\text{that is,}\qquad \alpha \leq \delta \leq \pi-\beta.
\]
If $\delta < \alpha$, or equivalently $\Delta > \cos(\alpha)$, $y_*$ is to the left of $y_h$ and the minimum of (\ref{eqn:simple})
in $[y_h,z_h]$ is attained at $y_h$. On the other hand, if $\delta > \pi-\beta$, or equivalently $\Delta < \cos(\pi-\beta)$,
$y_*$ is to the right of $z_h$ and the minimum of (\ref{eqn:simple})
in $[y_h,z_h]$ is attained at $z_h$.
\end{proof}

For the general case we simple apply the triangular update formula of Lemma~\ref{lem:update} to the
image of the triangle $\sigma_i$ under the linear transformation $M(x_h)^{-1/2}$.
This way we immediately obtain the following update procedure, writing $\langle p,q\rangle_{x} =
\langle p,q\rangle_{M(x)^{-1}}$, $\|p\|_{x} = \|p\|_{M(x)^{-1}}$,
$c_\alpha=\cos(\alpha)$, and $c_\beta=\cos(\beta)$ for short (note that we used the addition formula
to spell out $\cos(\alpha-\delta)$ for implementation purposes):

\begin{align*}
&\Delta = \dfrac{u_h(z_h)-u_h(y_h)}{\|z_h-y_h\|_{x_h}};\\*[2mm]
&c_\alpha = \dfrac{\langle x_h-y_h,z_h-y_h\rangle_{x_h}}{\|x_h-y_h\|_{x_h}\cdot\|z_h-y_h\|_{x_h}};\quad
c_\beta = \dfrac{\langle x_h-z_h,y_h-z_h\rangle_{x_h}}{\|x_h-z_h\|_{x_h}\cdot\|y_h-z_h\|_{x_h}};\\*[2mm]
&\texttt{if}\;\;  c_\alpha\leq \Delta \\*[1mm]
&\qquad u_i=u_h(y_h) + \|x_h-y_h\|_{x_h}; \\*[1mm]
&\texttt{else if}\;\;  \Delta\leq -c_\beta \\*[1mm]
&\qquad u_i=u_h(z_h) + \|x_h-z_h\|_{x_h}; \\*[1mm]
&\texttt{else} \\*[1mm]
&\qquad u_i=u_h(y_h) + \left(c_\alpha\Delta  + \sqrt{(1-c_\alpha^2)(1-\Delta^2)}\right) \|x_h-y_h\|_{x_h};
\end{align*}

\paragraph{\protect\emph{Remark}} With different ideas on a discretization, exactly the same update formula
has been obtained for the (classic) eikonal equation by
Kimmel and Sethian \cite{0908.65049} (see also \cite[{\S}10.3.1]{Sethian}), who use for acute triangulations the methodology
of \cite{BS} to construct upwind schemes on unstructured meshes, and, independently, by the geophysicist Fomel \cite{Fomel}, who
locally uses Fermat's principle of shortest travel times (which is closely related to our local use of the Hopf--Lax formula).

Sethian \cite[{\S}10.1]{Sethian} shows further that this update formula generalizes the upwind finite-difference scheme on structured grids
given by Rouy and Tourin \cite{0754.65069}.

\medskip\bigskip\section{Solving the Discrete System}\label{sect:gs}

\paragraph{\emph{A Review of Methods.}}

Theorem~\ref{thm:exiFEM} shows that the nonlinear discrete system (\ref{eq:fixedpoint}) can be solved
by the fixed-point iteration
$$
u_h^{n+1} = \Lambda_h u_h^n, \quad n=0,1,2,\ldots,
$$
for a suitably chosen initial iterate $u^0_h$. Such a fixed-point iteration uses the updated values at a nodal point $x_j$
only after all the updated values have been calculated. This corresponds to the classic {\em Jacobi--iteration} for linear systems
of equations and lends itself to direct parallelization.

If we sequentially traverse the nodal points in a given order and modify the iteration to always use the most recently
updated value, we obtain a nonlinear variant of the \emph{Gauss--Seidel iteration}. Rouy and Tourin \cite{0754.65069} used
such a nonlinear Gauss-Seidel iteration to solve a finite difference
discretization of the eikonal equation $\norm{Du(x)}=n(x)$ on a structured mesh.

For both iterative methods the complexity will typically scale as $O(N^{1+1/d})$,
where $N$ denotes the number of nodal points and $d$ the space dimension. This is because the information about the solution,
inherent initially to the boundary only, travels by next neighbor interaction at each run trough all nodal points. To spread
that information to the whole computational domain about $O(N^{1/d})$ runs are necessary. However, even though this heuristic
well explains the experimental observations, to our knowledge there is no rigorous proof of that in the literature.

In 1995 Sethian \cite{SethActa} and Tsitsiklis \cite{0831.93028} have shown independently that for eikonal equations on structured meshes
the nonlinear equation can be solved {\em exactly} in a \emph{single pass}, that is, by traversing the grid once using local
operations only. This {\em fast marching method} was later generalized to triangular meshes by Kimmel and Sethian \cite{0908.65049}.
It relies on the \emph{causality property}, namely that on an acute triangulation the value $u_h(x_h)$
depends only on the values in neighboring nodal points $y_h$ that are lower, $u_h(y_h)\leq u_h(x_h)$. So the discrete solution
can be computed starting from the neighborhood of the boundary moving further inwards the computational domain
along increasing values of $u_h$. However, on non-acute triangulations additional effort is necessary
to deal with the loss of this causality property (see \cite{0908.65049} for details).
The complexity of this method is $O(N\log(N))$ where the logarithm comes from administering a priority queue
of candidates for the next lowest value of $u_h$, such as a heap data structure.

For the particular Hamilton--Jacobi--Bellman equation (\ref{eq:HJB}) a single pass algorithm generalizing the fast marching method,
called the
\emph{ordered upwind method} (OUM), was introduced by Sethian and Vladimirsky \cite{0963.65076} and is
discussed in detail in \cite{SETHVLADSIAM}. This method is not a fast solver for a given discretization, but the
discretization is specifically designed for the needs of the fast solver.
Like the Hopf--Lax update the update formulas of the OUM are based on local variational principles (Bellman's principle).
However though, the OUM does
not solve (\ref{eq:fixedpoint}), since the update in $x_h$ is not necessarily computed from the neighborhood $\omega_h(x_h)$ but
from larger neighborhoods within the radius $h\cdot\nu$.
Here  $\nu=F^*/F_*$ denotes the \emph{anisotropy coefficient} of the Hamilton--Jacobi--Bellman equation (see (\ref{eq:Fstar})).
This quantity $\nu$ affects not only the complexity of the OUM, which is $O(\nu^{d-1} N \log(N))$,
but also its accuracy.

\smallskip
\paragraph{\emph{Adaptive Nonlinear Gauss--Seidel Iteration.}}
In this paper we propose an adaptive version of the nonlinear Gauss-Seidel iteration, which is a modelled after a
similar relaxation method \cite{Rude} for the multilevel solution of elliptic boundary value problems.
It turns out to be substantially faster than the standard Gauss-Seidel iteration, easy to implement and universal.

The adaptive Gauss-Seidel iteration differs from the standard one in two respects. First, like in the fast marching method,
only those nodal points are updated that ``have the
information'', that is, are neighbors of recently updated points. Second, the order of updates is not fixed but
varies as the iteration proceeds. Thus, a queue denoted by $\mathcal{Q}$ is administered
to provide the ordering of updates. However, other than in the fast-marching method where using the causality property
requires to
keep control of the point with minimal function value, the queue is now simply given the structure of a FIFO (first in first out)
stack: the nodal point staying longest in the queue is updated next.

The algorithm is passed a user-defined tolerance $\tol$ and it ends up with an approximate
finite-element solution $u_h \in V_h$ such that
$$
\norm{ u_h - \Lambda_h  u_h}_\infty \leq \tol.
$$
It is organized as follows:

\medskip

\begin{compactenum}
\item (Initialization) Let $u_h|_{\partial\Omega_h}=g|_{\partial\Omega_h}$, $u_h|_{\Omega_h}\equiv \infty$.\footnote{In fact any
value larger than the bound (\ref{eq:bound}) will do. However, taking $\infty$ makes the argument more elegant and can correctly
be implemented in IEEE arithmetic.}
Let $\mathcal{Q}$ be the list of all points $x_h\in\Omega_h$ that are adjacent to some boundary point
(in an arbitrary but fixed order).
\item (Iteration) Remove the first point $x_h$ from $\mathcal{Q}$ and compute the update value
\mbox{$u_{\rm new}=(\Lambda_h u_h)(x_h)$}.
\item If $ |u_{\rm new} - u_h(x_h)| > \tol$ then update $u_h(x_h)=u_{\rm new}$ and append all not yet enqueued
neighbors $y_h$ of $x_h$ to the queue $\mathcal{Q}$.
\item If $\mathcal{Q}\not = \emptyset$, goto (2).
\end{compactenum}

\medskip

To prove the convergence of this method, we denote the initial finite-element function of step 1 by $u_h^0$.
After the $n$th update has been performed in step 3 the actual finite-element function will be denoted by $u_h^n$.

\begin{theorem}
The algorithm generates a sequence $u_h^0,u_h^1,\ldots$, that is monotonously decreasing.
It terminates after finitely many steps with an approximate finite element solution  $u_h \in V_h$, such that
$\norm{u_h - \Lambda_h u_h}\leq \tol$.
\end{theorem}
\begin{proof}
The initialization $u_h|_{\Omega_h}\equiv \infty$ ensures that every point $x_h$ is updated at least once, as
the residual is $\infty$ when the first update value in $x_h$ is computed.
After the first update, $u_h(x_h)$ is assigned a finite value, since $x_h$ has a neighbor in $\partial\Omega_h$ or a neighbor, for which
a finitely valued update has already been computed. By induction on $n$ we get that at each later update of a nodal point $x_h$, all
neighbors of $x_h$ that have been changed over the last update can only have been assigned a lower value of $u_h$.
>From the monotonicity of $\Lambda_h$ we thus get the first assertion.

Since an update in step (3) only affects the residual in the neighboring points, which are
are immediately enqueued, it holds that
\[
\{ x_h \in\Omega_h \ : \ u_h^n(x_h)<\infty \textrm{~and~} |\Lambda_h u_h^n - u_h^n| > \tol \} \ \subset \ \mathcal{Q}
\]
for every $n\geq 0$. So if the algorithm terminates with $\mathcal{Q}=\emptyset$, the tolerance has been reached.

Otherwise, if the iteration does not terminate, then there is at least one nodal point $x_h^*$ that appears infinitely often
as the first element of the queue $\mathcal{Q}$ and gets updated at steps $n_j\to\infty$, $j \to \infty$. Hence,
there must be $|u_h^{n_j}(x_h^*)- u_h^{n_j-1}(x_h^*)|>\tol$ in contradiction to the convergence of $u_h^{n_j}(x_h^*)$ as $j\to\infty$ which
is implied by the monotonicity and the trivial lower bound $u_h^n \geq \min_{x\in\partial\Omega} g(x)$.
\end{proof}

Though the run-time complexity of the adaptive Gauss--Seidel iteration behaves probably at worst as $O(N^{1+1/d})$ like
in the standard Gauss-Seidel iteration, a lot of unnecessary updates are saved as we will see in the numerical
experiments of the next section.

\begin{figure}[tbp]
\includegraphics[width=0.39\textwidth]{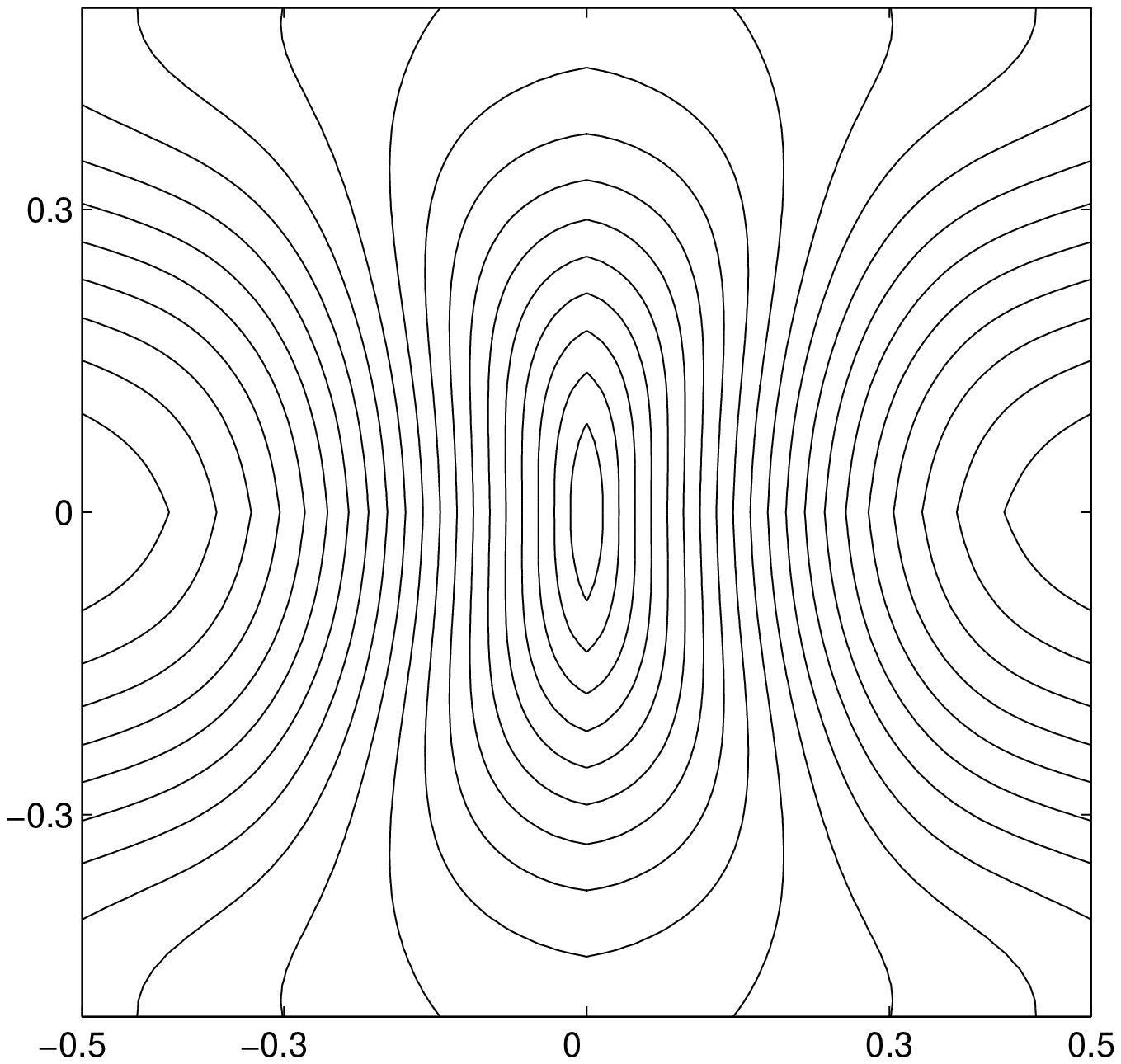}
\qquad
\includegraphics[width=0.46\textwidth]{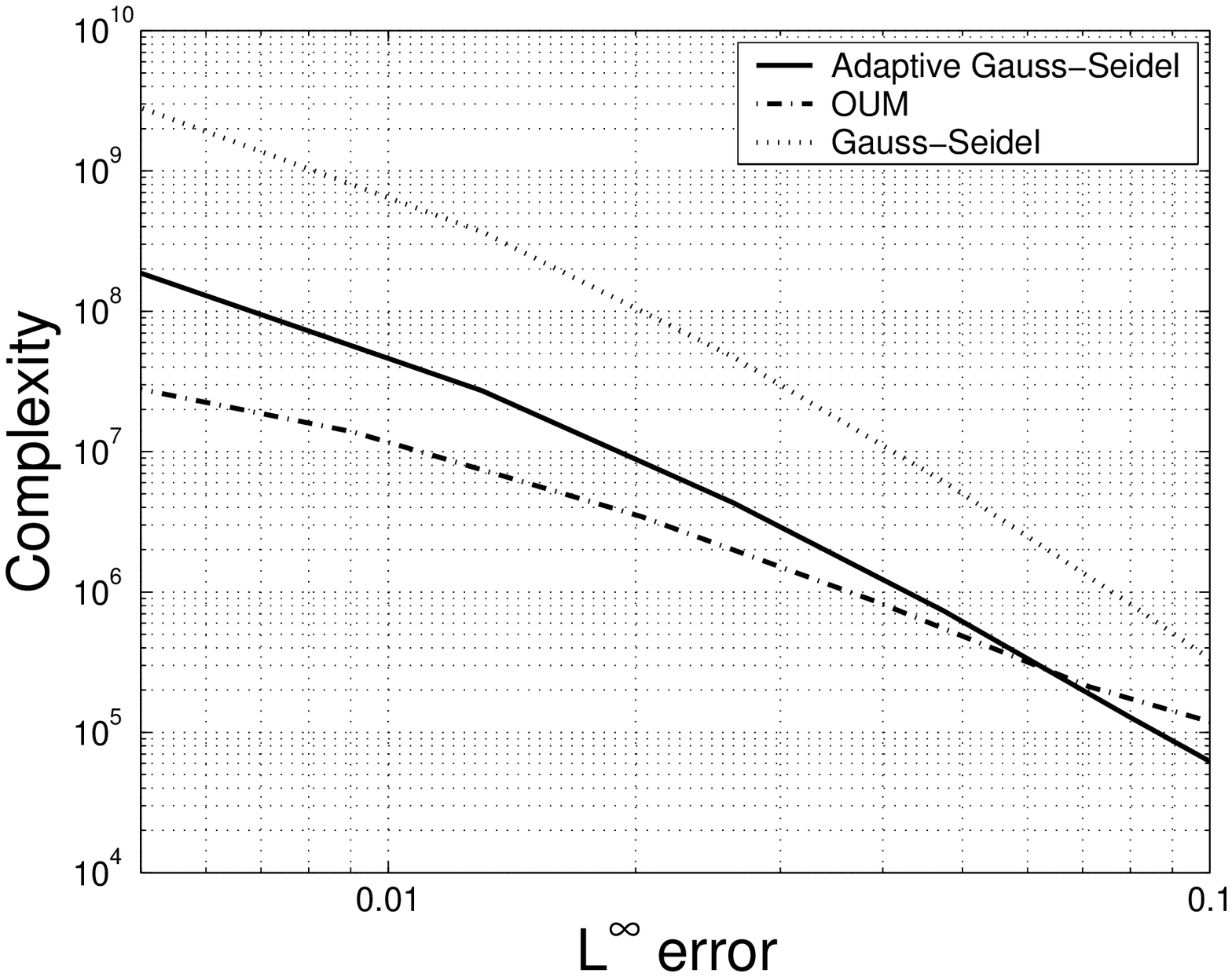}
\caption{{\protect\small Left: Contour plot of the distance function over the parameter plane.
Right: Accuracy/complexity of adaptive Gauss-Seidel iteration in comparison to the Gauss-Seidel iteration and OUM.}}
\label{fig:FIG2}
\end{figure}

\medskip\bigskip\section{Numerical Experiments}\label{sect:exp}

For the two following examples the solutions were computed on unstructured meshes with $23^2$, $45^2$, $91^2$, $181^2$, and $725^2$ nodal points.
A solution on a mesh with $1451^2$ points served to estimate the discretization error. The iterative methods were used with an absolute
tolerance $\tol = 10^{-8}$.

The first example concerns the distance map on the torus given by the immersion
\[
f(x_1,x_2)=  \left(\cos(2 \pi x_1)(5+4\cos(2 \pi x_2)),
\sin(2\pi x_1)(5+4\cos(2\pi x_2)), \sin(2\pi x_2) \right).
\]
Using the Gram matrix $G(x)=Df(x)^{T}Df(x)$ and $\rho^2(x,q)=\ska{q}{G(x)q}$
the distance
between the points $f(x)$ and $f(y)$ on the manifold is given by the function $\delta(x,y)$ as defined in (\ref{eqn:DELTA}).
With the results of {\S}\ref{sect:eikonal} and \cite[Thm.~5.3(iv)]{0497.35001}
we obtain that $u(x)=\delta(x,0)$ is the viscosity solution of the Dirichlet problem
$$
\norm{Du}_{G(x)^{-1}}=1 \textrm{~on~}\Omega\setminus\{0\}, \qquad u(0)=0,
$$
where $\Omega=[-0.5,0.5]^2$. The solution is shown to the left of Figure~\ref{fig:FIG2} as a contour plot.

To the right of  Figure~\ref{fig:FIG2} we compare the accuracy and complexity of the adaptive Gauss-Seidel method with both the standard
Gauss-Seidel iteration and the OUM.\footnote{We have coded the OUM from \cite{SETHVLADSIAM}
with a little completion that turned out to be necessary:
Considered points have also to be updated, if they depend on an edge that drops out of
the accepted front. If some point $x_h$ gets accepted this may happen to any edge opposite to $x_h$ in
$\omega_h(x_h)$.} Here, by \emph{complexity} we mean the total
number of updates calculated on a triangle by a formula such as the one at the end of {\S}\ref{sect:eikonal}.
We observe that the adaptive Gauss--Seidel iteration is more than a factor of 10 faster than the standard
Gauss--Seidel iteration but displays the same asymptotic rate of complexity. The OUM show, as theoretically expected,
a better rate of complexity that, in this example, becomes significant even at larger tolerances.

\medskip

\begin{figure}
\includegraphics[width=0.39\textwidth]{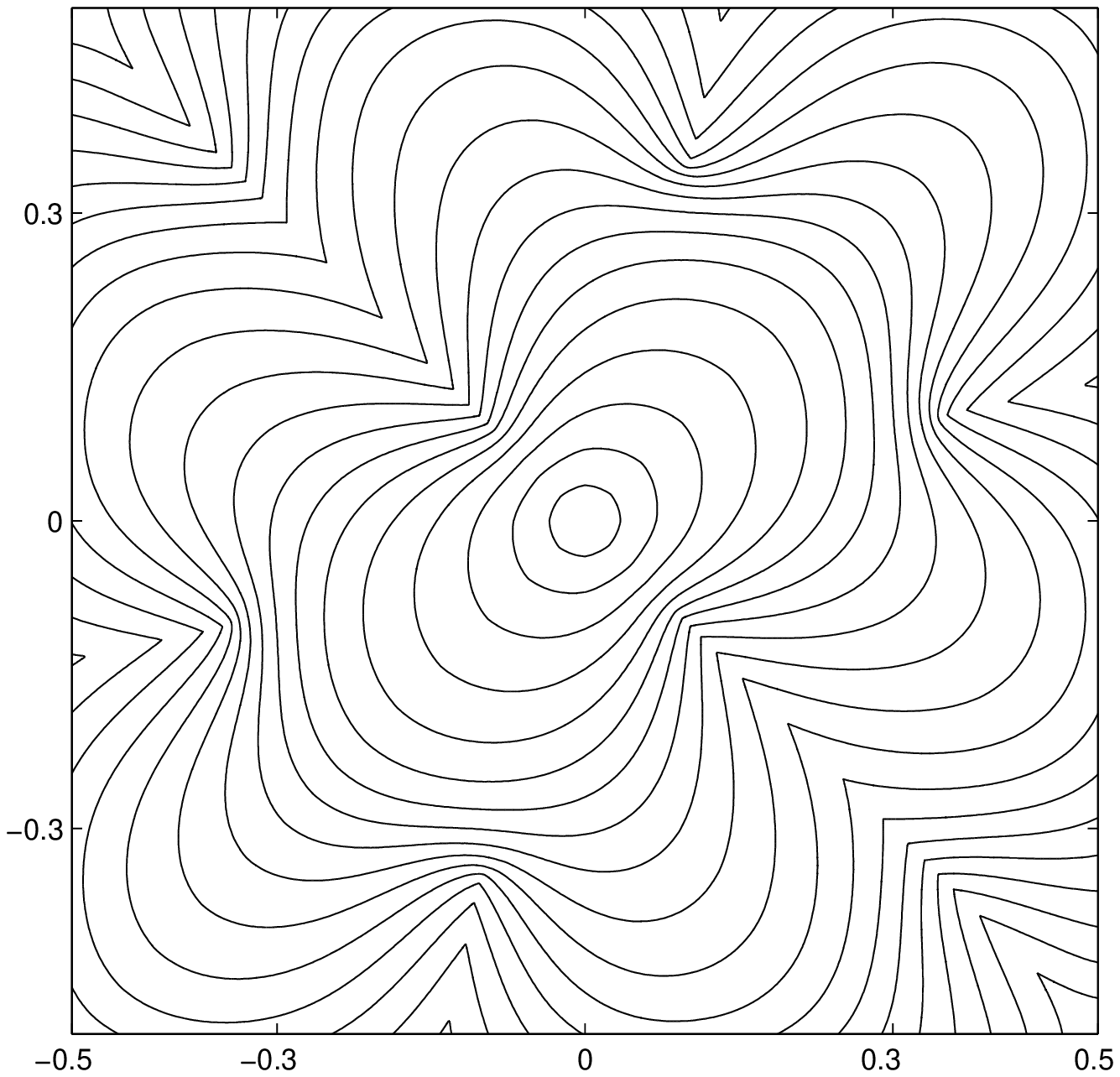}
\qquad
\includegraphics[width=0.46\textwidth]{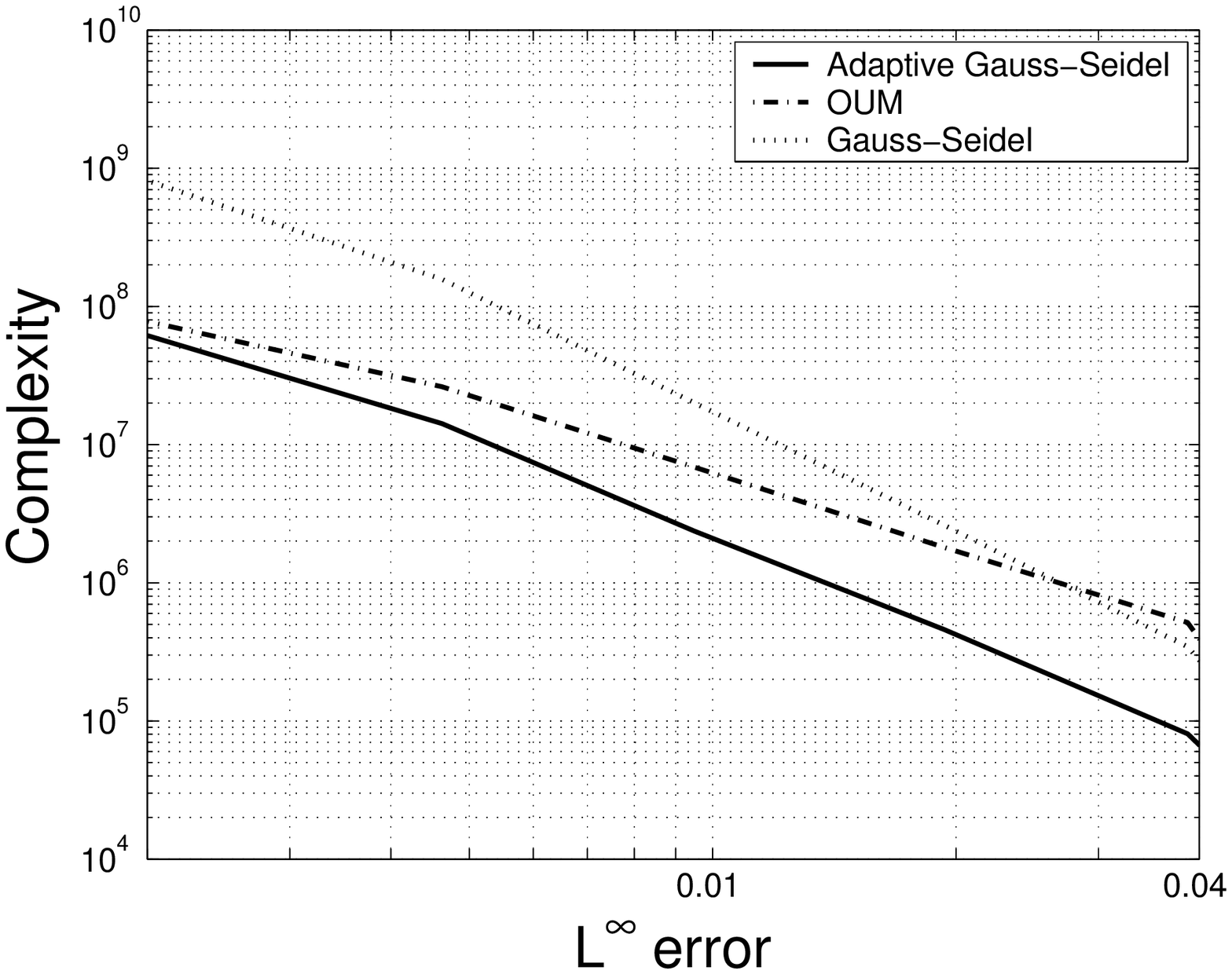}
\caption{{\protect\small Left: Value function of some min-time optimal control problem. Right: Complexity/accuracy
of the methods in comparison.}}
\label{fig:fig4}
\end{figure}

The second example is taken from \cite{1002.65112} and shows the effect of a moderately large anisotropy coefficient
$\nu=19$, which affects the complexity of the OUM.
We consider a simple min-time optimal control problem
governed by the dynamical system
\[
y'(t)=a(t) + b(y(t)), \qquad y(0)=x.
\]
The controls $a(\cdot)$ are taken from $\mathcal{A}=\{ a:[0,\infty) \rightarrow S^1 \textrm{~measurable} \}$
and
\[
b(y)=-0.9\sin(4\pi y_1) \sin(4\pi y_2) \cdot \frac{y}{\norm{y}}.
\]
For $x\in\Omega=[-0.5,0.5]^2$ and a control $a$ we denote by $T_x(a)$ the minimal time that
the trajectory $y(\cdot)$ takes to reach the origin.

Following \cite[p.~241, Thm.~2.6]{0890.49011} the value function $u(x)=\inf_{a\in\mathcal{A}} T_x(a)$ is the
viscosity solution of the Hamilton--Jacobi--Bellman equation
\[
H(x,Du)=\max_{\|a\|=1}  \ska{a + b(x)}{-Du} - 1  = 0 , \qquad u(0)=0.
\]
One figures out that $H(x,p)=\norm{p}-\ska{b(x)}{p}-1$; (H1), (H2), (H4$'$), and, as $\norm{b}\leq 0.9$, the
coercivity condition (H3) are fullfilled. A short calculation shows that
\[
\rho(x,q)=\frac{\norm{q}}{\left( 1-\norm{b(x)}^2 + \ska{b(x)}{q/\|q\|}^2\right)^{1/2}-
\ska{b(x)}{q/\|q\|}},
\]
compare also (\ref{eq:rhoHJB}) and \cite[Eq.~(19)]{1002.65112}.
The solution calculated on a $253\times 253$ mesh can be found to the left of Figure~\ref{fig:fig4}. To the right of this figure
the accuracy of the approximate finite-element solution is shown versus the complexity of the iteration.
A comparison of the (adaptive) Gauss--Seidel iteration with the OUM is shown to the right of Figure~\ref{fig:fig4}.
Again we observe that the adaptive variant of the Gauss--Seidel iteration is by about a factor of ten more efficient
than the standard one. This time, however, the quite sophisticated order upwind method
behaves less favorable: because of the large anisotropy coefficient the break-even point at which the OUM
becomes more efficient than the simple adaptive Gauss--Seidel
iteration is at a mesh-size of more than $725^2 = 525,\!625$ nodal points. We expect this effect to become even more pronounced in 3D and higher,
because of the increasingly better complexity rate of the Gauss--Seidel iteration.

\providecommand{\bysame}{\leavevmode\hbox to3em{\hrulefill}\thinspace}
\providecommand{\MR}[1]{\relax\ifhmode\unskip\space\fi MR #1.}
\providecommand{\MRhref}[2]{%
  \href{http://www.ams.org/mathscinet-getitem?mr=#1}{#2}
}
\providecommand{\href}[2]{#2}


\begin{thebibliography}{BCD97}

\bibitem[Alt99]{Alt02}
Hans~Wilhelm Alt, \emph{Lineare {F}unktionalanalysis}, third ed.,
  Springer-Verlag, Berlin, 1999.

\bibitem[BCD97]{0890.49011}
Martino Bardi and Italo Capuzzo-Dolcetta, \emph{Optimal control and viscosity
  solutions of {H}amilton-{J}acobi-{B}ellman equations}, Birkh{\"a}user,
  Boston, 1997. \MR{99e:49001}

\bibitem[BS98]{BS}
Timothy~J. Barth and James~A. Sethian, \emph{Numerical schemes for the
  {H}amilton-{J}acobi and level set equations on triangulated domains}, J.
  Comput. Phys. \textbf{145} (1998), no.~1, 1--40. \MR{99d:65277}

\bibitem[CEL84]{0543.35011}
Michael~G. Crandall, Lawrence~C. Evans, and Pierre-Louis Lions, \emph{Some
  properties of viscosity solutions of {H}amilton-{J}acobi equations}, Trans.
  Amer. Math. Soc. \textbf{282} (1984), no.~2, 487--502. \MR{86a:35031}

\bibitem[Eva98]{Evans}
Lawrence~C. Evans, \emph{Partial differential equations}, American Mathematical
  Society, Providence, 1998. \MR{99e:35001}

\bibitem[Fom97]{Fomel}
Sergey Fomel, \emph{A variational formulation of the fast marching eikonal
  solver}, Tech. Report 95, pp.~127--149, Stanford Exploration Project,
  Stanford University, 1997,
  \\\href{http://sepwww.stanford.edu/public/docs/sep95/toc_html/}{sepwww.stanf%
ord.edu/public/docs/}.

\bibitem[Ish87]{0644.35017}
Hitoshi Ishii, \emph{A simple, direct proof of uniqueness for solutions of the
  {H}amilton-{J}acobi equations of eikonal type}, Proc. Amer. Math. Soc.
  \textbf{100} (1987), no.~2, 247--251. \MR{88d:35040}

\bibitem[KS98]{0908.65049}
Ron Kimmel and James~A. Sethian, \emph{Computing geodesic paths on manifolds},
  Proc. Natl. Acad. Sci. USA \textbf{95} (1998), no.~15, 8431--8435.
  \MR{99d:65359}

\bibitem[Lio82]{0497.35001}
Pierre-Louis Lions, \emph{Generalized solutions of {H}amilton-{J}acobi
  equations}, Pitman, Boston, 1982. \MR{84a:49038}

\bibitem[LYC03]{Xiang}
Xiang-Gui Li, Wei Yan, and C.~K. Chan, \emph{Numerical schemes for
  {H}amilton-{J}acobi equations on unstructured meshes}, Numer. Math.
  \textbf{94} (2003), no.~2, 315--331. \MR{2004b:65153}

\bibitem[OF03]{Osher}
Stanley Osher and Ronald Fedkiw, \emph{Level set methods and dynamic implicit
  surfaces}, Springer-Verlag, New York, 2003. \MR{2003j:65002}

\bibitem[OS88]{0659.65132}
Stanley Osher and James~A. Sethian, \emph{Fronts propagating with
  curvature-dependent speed: algorithms based on {H}amilton-{J}acobi
  formulations}, J. Comput. Phys. \textbf{79} (1988), no.~1, 12--49.
  \MR{89h:80012}

\bibitem[PR93]{Rude}
Christoph Plaum and Ulrich R{\"u}de, \emph{Gau{\ss}' adaptive relaxation for
  the multilevel solution of partial differential equations on sparse grids},
  Tech. Report SFB-Bericht 342/13/93, Technische Universit{\"a}t M{\"u}nchen,
  1993,
  \href{http://www10.informatik.uni-erlangen.de/~ruede/Preprints/preprints.htm%
l}{www10.informatik.uni-erlangen.de/\~\/ruede/}.

\bibitem[Roc70]{0932.90001}
R.~Tyrrell Rockafellar, \emph{Convex analysis}, Princeton Univ. Press,
  Princeton, 1970. \MR{43:445}

\bibitem[RT92]{0754.65069}
Elisabeth Rouy and Agn{\`e}s Tourin, \emph{A viscosity solutions approach to
  shape-from-shading}, SIAM J. Numer. Anal. \textbf{29} (1992), no.~3,
  867--884. \MR{93d:65019}

\bibitem[Set96]{SethActa}
James~A. Sethian, \emph{Theory, algorithms, and applications of level set
  methods for propagating interfaces}, Acta numerica, 1996, Acta Numer.,
  vol.~5, Cambridge Univ. Press, Cambridge, 1996, pp.~309--395. \MR{99d:65397}

\bibitem[Set99]{Sethian}
\bysame, \emph{Level set methods and fast marching methods}, second ed.,
  Cambridge Monographs on Applied and Computational Mathematics, vol.~3,
  Cambridge University Press, Cambridge, 1999, Evolving interfaces in
  computational geometry, fluid mechanics, computer vision, and materials
  science. \MR{2000c:65015}

\bibitem[SV00]{0963.65076}
James~A. Sethian and Alexander Vladimirsky, \emph{Fast methods for the eikonal
  and related {H}amilton-{J}acobi equations on unstructured meshes}, Proc.
  Natl. Acad. Sci. USA \textbf{97} (2000), no.~11, 5699--5703. \MR{2001b:65100}

\bibitem[SV01]{1002.65112}
\bysame, \emph{Ordered upwind methods for static {H}amilton-{J}acobi
  equations}, Proc. Natl. Acad. Sci. USA \textbf{98} (2001), no.~20,
  11069--11074. \MR{2002g:65133}

\bibitem[SV03]{SETHVLADSIAM}
\bysame, \emph{Ordered upwind methods for static {H}amilton-{J}acobi equations:
  theory and algorithms}, SIAM J. Numer. Anal. \textbf{41} (2003), no.~1,
  325--363. \MR{1 974 505}

\bibitem[Tsi95]{0831.93028}
John~N. Tsitsiklis, \emph{Efficient algorithms for globally optimal
  trajectories}, IEEE Trans. Automat. Control \textbf{40} (1995), no.~9,
  1528--1538. \MR{96d:49039}

\end{thebibliography}
\end{document}